\documentclass[11pt]{amsart}
\usepackage{amssymb, amscd}
\usepackage{latexsym,epsfig}
\usepackage[all]{xy}
\usepackage{pst-all}
\numberwithin{equation}{section}

\newcommand{\VVl}{\VV_{\lambda}}
\newcommand{\Pa}{\lambda = \{\lambda_1 \geq \lambda_2 \geq \lambda_3\geq 0\}}
\newcommand{\lam}{\lambda}
\newcommand\la[1]{\lambda_{{#1}}} 
\newcommand{\Sy}{{\mathrm{Sym}}} \newcommand{\we}{\wedge}
\newcommand{\noir}{\Sy^{\la1-\la2}\VV\otimes
\Sy^{\la2-\la3}\left(\we^2\VV\right) \otimes
\Sy^{\la3}\left(\we^3\VV\right)}
\newcommand{\VV}{\mathbb V}

\newcommand{\FF}{\mathbb F}
\newcommand{\eul}{e_c(\md,\VVl^{\prime})}

 \newcommand{\hy}{{\mathcal H}_3}
\newcommand{\md}{{\mathcal M}_3} \newcommand{\mde}{{\mathcal M}_{3,1}}

 \newcommand{\CC}{\mathbb C}

\newcommand{\QQ}{\mathbb Q}
\newcommand{\ZZ}{\mathbb Z}

\newcommand{\Sp}{{\mathrm{Sp}}(6,\QQ)}

\newcommand{\A}{\mathcal A}

\newtheorem{theorem}{Theorem}[section]
\newtheorem{lemma}[theorem]{Lemma}
\newtheorem{proposition}[theorem]{Proposition}
\newtheorem{corollary}[theorem]{Corollary}

\newtheorem{definition-lemma}[theorem]{Definition-Lemma}

\theoremstyle{definition}

\theoremstyle{remark}
\newtheorem{remark}[theorem]{Remark}

\begin{document}

\title[]{The Euler Characteristic of local systems on 
\\ the moduli of curves and abelian varieties of genus three} 
\author{Jonas Bergstr\"om}
\address{Institut Mittag-Leffler, Aurav\"agen 17, 18260 Djursholm, Sweden}
\email{jonasb@math.kth.se}

\author{Gerard van der Geer}
\address{Korteweg-de Vries Instituut, Universiteit van
Amsterdam, Plantage \newline Muidergracht 24, 1018 TV Amsterdam, The Netherlands.}
\email{geer@science.uva.nl}


\begin{abstract}
We show how to calculate the Euler characteristic of a local system $\VVl$ 
associated to an irreducible representation $V_{\lambda}$ of the symplectic
group of genus $3$ on the moduli space $\md$ of curves of genus $3$
and the moduli space $\A_3$ of principally polarized 
abelian varieties of dimension~$3$.
\end{abstract}

\maketitle

\begin{section}{Introduction}
\label{sec: intro}
An irreducible representation $V_{\lambda}$ with highest weight
$\lambda$ of the symplectic group ${\rm GSp}(2g,\QQ)$ defines a
local system $\VVl$ on the moduli space $\A_g$ of principally
polarized abelian varieties. Pull back under the Torelli morphism
defines a local system  $\VVl^{\prime}$ on the moduli space ${\mathcal M}_g$
of curves of genus $g$. The cohomology of such local systems is
intimately connected with the cohomology of the moduli spaces
${\mathcal M}_{g,n}$ of $n$-pointed curves and also with vector-valued
Siegel modular forms of degree $g$, cf., \cite{del, G2, FvdG}. 
Therefore it is of some interest
to be able to calculate the Euler characteristic of such a
local system. 

Getzler showed in \cite{G} how to do this for 
${\mathcal M}_2$ and Bini and the second author did this for 
local systems on the hyperelliptic locus $\hy$ 
of genus $3$ in \cite{BvdG} by calculating the Euler numbers
of a stratification.
In the case of $\md^0$, the non-hyperelliptic locus of $\md$,
it seems difficult to calculate the Euler characteristics of the
individual strata directly. Instead, we use the 
results of the first author (see \cite{B1})
on the motivic Euler characteristics of $\VVl^{\prime}$ 
on $\md$ for $\lambda$ of small weight, obtained by counting points
over finite fields, to calculate linear relations
between the Euler numbers of the strata. 
We were not able to determine the Euler numbers of all
strata, but the information obtained suffices for determining
the Euler characteristic of all local systems on $\md$. To go from 
$\md$ to $\A_3$ we need to calculate the Euler characteristics
of the pull backs of our local systems to ${\mathcal M}_2 \times 
\A_1$ and to the locus $\A_{1,1,1}$
of products of three elliptic curves.  
We illustrate our results by giving tables of Euler characteristics
on $\md$ and $\A_3$.
We hope that this is a step towards a better
understanding of the cohomology of local systems of genus~$3$.

We thank Torsten Ekedahl for providing us with the
argument for the invariance of the Euler characteristic in Section 
\ref{basechange}.
We also thank the Mittag-Leffler Institute for the hospitality enjoyed 
during the preparation of this paper.

\end{section}
\begin{section}{The Euler characteristic of a local system}
\label{section1}
Let $\md$ be the moduli space of smooth genus $3$ curves. It is
a Deligne-Mumford stack of relative dimension $6$ over $\ZZ$. 
The universal curve $\pi' : \mde \to \md$ defines,
for any prime $\ell$, a natural $\ell$-adic local system 
$R^1\pi'_*\QQ_{\ell}$ of rank $6$ on $\md \otimes \ZZ[1/\ell]$.
We shall denote this local system by 
${}_{\ell}\VV^{\prime}$, or simply by $\VV^{\prime}$. 
It comes equipped with a non-degenerate symplectic pairing 
$\VV^{\prime} \times \VV^{\prime} \to \QQ_{\ell}(-1)$.

Similarly, we have the moduli space $\A_3$
of principally polarized abelian varieties
of dimension $3$, again a Deligne-Mumford stack of relative dimension $6$
over $\ZZ$. The universal abelian threefold 
$\pi\colon{\mathcal X}_3 \to \A_3$ also defines a natural $\ell$-adic local system
$R^1\pi_* \QQ_{\ell}$ on $\A_3 \otimes \ZZ[1/\ell]$,
which we denote by ${}_{\ell} \VV$, or simply by $\VV$.
There is the Torelli morphism $t_3\colon \md \to \A_3$
of degree $2$ between the stacks.
The local system $\VV^{\prime}$ on $\md\otimes \ZZ[1/\ell]$ 
is a pull back from the local system $\VV$ on $\A_3\otimes \ZZ[1/\ell]$.

A partition $\Pa$ of weight $|\lambda|=\lambda_1+
\lambda_2+\lambda_3$ determines an irreducible representation 
$V_{\lambda}$ 
of $\Sp$ associated to $\lam$. We lift it to a representation 
of ${\rm GSp}(6,\QQ)$ with dominant weight
$(\lambda_1-\lambda_2)\gamma_a+(\lambda_1-\lambda_2)\gamma_b
+\lambda_3 \gamma_c-|\lambda| \eta$ with $\gamma_a$, $\gamma_b$
and $\gamma_c$ suitable fundamental roots and 
$\eta$ the multiplier representation.
Any such representation yields a symplectic local system 
$\VVl$ of weight $|\lambda|$ on $\A_3\otimes \ZZ[1/\ell]$, 
which appears `for the first time' in the decomposition of
$$
\noir,
$$
into irreducibles. 
If, for example, $\lambda= \{\la1 \geq 0\geq 0\}$, then
$\VV_{\lambda}=\Sy^{\la1}(\VV)$.
Similarly, we have a local system $\VVl^{\prime}$ on $\md \otimes \ZZ[1/\ell]$.

Our goal is to determine in an algorithmic way,
for $k=\CC$ or $k=\bar{\FF}_p$, 
the Euler characteristic of compactly 
supported cohomology of $\VVl^{\prime}$ on $\md$ and of $\VVl$ on $\A_3$.

Compactly supported cohomology of a local system on a 
Deligne-Mumford stack presents
various subtleties if the order of the automorphism groups is not
invertible on the base, cf., the discussion in \cite{L-MB},
Remarque 18.3.3. 
When the order of the automorphism groups is invertible one can 
manage with an \emph{ad hoc} approach and that is what we do here.

Suppose ${\mathcal X}$ is a Deligne-Mumford stack, 
$\pi\colon {\mathcal X} \to X$ the natural map to its coarse moduli space, 
$f\colon {\mathcal X} \to S$ a map to a scheme $S$
and $g\colon X \to S$ the factorisation $f=g\pi$ through the coarse
moduli space. For ${\mathcal F}$ a $\QQ_{\ell}$-sheaf (with $\ell$
invertible on our base) on ${\mathcal X}$ we have 
$R\pi_*{\mathcal F} = \pi_*{\mathcal F}$, since $\pi$ is finite,
and hence by the Leray spectral 
sequence we have $Rf_*{\mathcal F} = Rg_*(\pi_*{\mathcal F})$.  We define
the direct image with compact support by
$Rf_!{\mathcal F}:=Rg_!(\pi_*{\mathcal F}).$

With this definition we have Poincar\'e duality for a local $\QQ_\ell$-system 
$\mathcal F$ on $\mathcal X$, if $\mathcal X$ is assumed to be smooth over some 
spatial base $S$ and purely $d$-dimensional. 
That is, we have a natural isomorphism
$$R{\mathcal H}om_S(Rf_!\mathcal{F},\QQ_\ell) \cong 
Rf_*(\mathcal{F}^{\vee}(-d)[-2d]).$$ 
Indeed, by the duality theorem for $g$ we have
$$
R{\mathcal H}om_S(Rf_!\mathcal{F},\QQ_\ell) \cong 
Rg_*\bigl(R{\mathcal H}om_X(\pi_*\mathcal{F},g^!\QQ_\ell)\bigr).
$$
Now, $g \colon X \to S$ is locally (in the \'etale topology) a quotient $U/G \to S$ of a
smooth morphism $h\colon U \to S$ of schemes by the action of a finite group
$G$. Since $h$ is smooth we have $h^!\QQ_\ell \cong \QQ_\ell(-d)[-2d]$ and 
pushing down this sheaf to $U/G$ and
taking $G$-invariants gives $g^!\QQ_\ell \cong \QQ_\ell(-d)[-2d]$. 
Furthermore,
we may choose $U$ to be an \'etale neighbourhood of an arbitrary point of
$\mathcal X$ and since $RHom_{U/G}(\pi_*{\mathcal
F}|_{U/G},\QQ_\ell) \cong RHom_U(\pi_*{\mathcal F}|_U,\QQ_\ell)^G$ 
we see that $R{\mathcal
H}om_X(\pi_*{\mathcal F},\QQ_\ell) \cong \pi_* ({\mathcal F})^{\vee}$. 
But since $\pi_*{\mathcal F}|_U$ is a local system we find by taking
$G$-invariants that $\pi_* ({\mathcal F})^{\vee}=\pi_*({\mathcal F}^{\vee})$.
Putting these statements 
together gives us the desired duality formula. (Note that the sheaf
$\pi_*{\mathcal F}$ need not be a local system on $X$.)

We apply this to the stack $\md$ by taking the compactly supported cohomology of
the direct image $\nu_* \VVl^{\prime}$ on the coarse moduli space $M_3$ under
the natural map $\nu\colon \md \to M_3$ and similarly for $\A_3$. We thus write
$$
\label{eulcar}
e_c(\md \otimes k, \VVl^{\prime})=\sum_{i=0}^{12}(-1)^i 
\dim_{\QQ_{\ell}} H^i_c(\md \otimes k, \VVl^{\prime})
$$ 
and
$$
e_c(\A_3\otimes k,\VVl)=\sum_{i=0}^{12} (-1)^i \dim_{\QQ_{\ell}} H^i_c(\A_3\otimes k,\VVl),
$$
where $H^{\bullet}_c$ means compactly supported $\ell$-adic \'etale 
cohomology, and $\VVl$ means the local system ${}_{\ell} \VVl$ 
with $\ell$ a prime different from the characteristic of our field.
We shall see in the next paragraph that these Euler characteristics are 
independent of the choice of $\ell$, and in the next section
that they are independent of $k$, which justifies our notation.

In characteristic $0$ we also have the local system 
$R^1\pi'_*\QQ$ on $\md \otimes \QQ$ (resp. $R^1\pi_*\QQ$ on 
$\A_3 \otimes \QQ$), denoted ${}_0 \VV'$ (resp. ${}_0 \VV$) and for each 
$\lambda$, in the same 
way as in the $\ell$-adic case, an associated local system ${}_0 \VVl'$ 
(resp. ${}_0 \VVl$).
The comparison theorem 4.1 in \cite{SGA4XVI} tells us that the Euler 
characteristic 
$e_c(\md \otimes \CC, {}_0 \VVl')$ (resp. $e_c(\A_3\otimes \CC,{}_0 \VVl)$) 
with respect to compactly supported topological 
cohomology is equal to the $\ell$-adic versions 
$e_c(\md \otimes \CC, \VVl^{\prime})$ 
(resp. $e_c(\A_3\otimes \CC,\VVl)$) defined above.

\end{section}

\begin{section}{Base change}\label{basechange}
In this section we show that the $l$-adic variants of our local systems
commute with base change and we conclude that the Euler characteristic
$e_c(\md\otimes k, \VVl')$ is independent of $k$ being equal to
$\CC$ or $\bar{\FF}_p$, for any $p$ different from $\ell$.

\begin{theorem}\label{locallyconstant}
Let $f\colon X \to S$  be a smooth, proper, relatively Deligne-Mumford (cf.,
\cite[7.3.3]{L-MB}) map between algebraic stacks such that the prime
$\ell$ is invertible on $S$. Suppose given a
relative normal crossing divisor $D \subset X$ and a locally constant sheaf
${\mathcal E}$ of $\ZZ_{\ell}$, $\QQ_{\ell}$ or $\ell$-torsion modules
on $U:=X\setminus D$. If ${\mathcal E}$ is tamely 
ramified along $D$, then the 
$R^ig_*{\mathcal E}$ are locally constant and commute with base change where $g$
is the structure map $g\colon U \to S$.
\begin{proof}
We start by showing that if $j\colon U \to X$ is the open embedding, then
$Rj_*{\mathcal E}$ commutes with base change on $S$ and is
constructible. Indeed, this is a local statement on $X$ so we may assume that
$X$ is a scheme and then it follows from 
\cite[Th.\ finitude, 1.3.3]{SGA41/2}. For $f$ we 
have the proper base change theorem \cite[Thm.~18.5.1]{L-MB}, using
\cite[Thm.~16.6]{L-MB}. Hence, we conclude that $Rg_*{\mathcal
E} = Rf_*Rj_*{\mathcal E}$ commutes with base change and is constructible.

Furthermore, the smooth base change theorem \cite[Cor. 1.2]{SGA4XVI} for
$f$ is true as it is local on $X$. 
The fact that $R^ig_*{\mathcal E}$ is locally constant
now follows as in the proof of \cite[Thm 2.1]{SGA4XVI}. Indeed, properness is
used only through the proper base change theorem (in fact only to reduce to the
case when the base is normal) and, in turn, it is only applied to show that
$R^ig_*{\mathcal E}$ commutes with base change.
\end{proof}
\end{theorem}
\begin{corollary}\label{cor-locallyconstant}
Let $\rho\colon {\mathcal M}_g \to {\rm Spec}(\ZZ[1/\ell])$ be the structure map
and $\pi : {\mathcal M}_{g,1}\to {\mathcal M}_g$ the universal curve.  If
${\mathcal V}$ is a locally constant system on ${\mathcal M}_g$ then
$R^i\rho_!{\mathcal V}$ is a locally constant system over ${\rm
Spec}(\ZZ[1/\ell])$ that commutes with base change.
\begin{proof}
The statement follows from the corresponding statement for $R^i\rho_*{\mathcal
V}$ by duality (see Section 2) and the case of $R^i\rho_*{\mathcal V}$ follows
from Theorem~\ref{locallyconstant} applied to $\overline{\mathcal M}_g \to {\rm
Spec}(\ZZ[1/\ell])$ once one has verified tameness. By purity and Abhyankar's
lemma \cite[X.3.6]{SGA1}, it is enough to verify this at the localization of
$\overline{\mathcal M}_g$ at a generic point of a component of the
boundary. That generic point is however of characteristic~$0$.
\end{proof}
\end{corollary}
In view of this independence of the characteristic we shall often just
write $e_c(\md , \VVl^{\prime})$ instead of 
$e_c(\md \otimes k, \VVl^{\prime})$ 
and $e_c(\A_3,\VVl)$ instead of $e_c(\A_3\otimes k,\VVl)$.
The invariance of the Euler characteristic for $\A_3$ follows from the
stratification of $\A_3$ ($= t_3(\md) \, \sqcup \, t_2(\mathcal{M}_2) \times \A_1 \, \sqcup \, \A_{1,1,1}$) found in the beginning of Section \ref{eulA3}.

\section{Calculating the Euler characteristics} \label{calc}
We now turn to the calculation of the Euler characteristics. 
The first remark is that the Euler characteristics are unchanged 
if we change the representation $\lambda$ by a power of the multiplier 
$\eta$ (twisting). We can therefore restrict our attention from 
${\rm GSp}(6,\QQ)$ to ${\rm Sp}(6,\QQ)$.

In this section we will work over $\CC$ and write $\md$ for $\md \otimes \CC$. 
The local system $\VVl'$ will stand for ${}_0 \VVl'$.
The Euler characteristic $\eul$ is calculated by descending from 
the stack $\md$ to the coarse moduli space $M_3$ under the natural 
map $\mu\colon \md \to M_3$ and by calculating
$e_c(M_3,\mu_*(\VVl^{\prime}))$ using the stratification of $M_3\otimes \CC$ by the
automorphism group of the genus $3$ curve. 
On each stratum $\Sigma(G)$ the automorphism group 
of every curve is equal to $G$, and 
the direct image $\mu_*(\VV^{\prime})$ is a local system. The
Euler characteristic is then given by (see, e.g., \cite[Theorem 5.13]{Dav},
and also \cite{Illusie})
\begin{equation} \label{eq-eulstrat}
\eul = \sum_G e_c(\Sigma(G))\dim({\VVl^{\prime}}^G),
\end{equation}
where $e_c(\Sigma(G))$ is the topological Euler characteristic of
$\Sigma(G)$ and ${\VVl^{\prime}}^G$ is the space of $G$-invariants. 

The dimension of the invariant part 
${\VVl^{\prime}}^G$ is determined as follows.
The action of $G$ on the cohomology group $H^1(C,\QQ)$ of a curve $C$
of genus $3$ defines a homomorphism
$r \colon G \to {\rm Sp}(6,\QQ)$. We denote the eigenvalues of
$r(\gamma)$ by 
$$
\{a_{\gamma},b_{\gamma},c_{\gamma},a_{\gamma}^{-1},b_{\gamma}^{-1},
c_{\gamma}^{-1}\}.
$$
Let $h_d$ be the complete symmetric function of degree $d$ in six
variables. For the proof of the following fact we refer to
\cite[Prop.\ 24.22]{FH}.

\begin{proposition}\label{diminv}
If $J_{\lambda}$ is the determinant of the $3 \times 3$ matrix
whose $i$-th row is 
$$
(J_{\lambda_i-i+2}, \, 
J_{\lambda_i-i+2}+J_{\lambda_i-i},
\, J_{\lambda_i-i-1})
$$
and 
$$
J_d(a_{\gamma},b_{\gamma},c_{\gamma}):=h_d(a_{\gamma},b_{\gamma},c_{\gamma},a_{\gamma}^{-1},b_{\gamma}^{-1},
c_{\gamma}^{-1}),
$$
then we have
$$
\dim ({\VVl^{\prime}}^G)= \frac{1}{\# G} \sum_{\gamma \in G}
J_{\lambda}(a_{\gamma},b_{\gamma},c_{\gamma}).
$$
\end{proposition}
Hence, we can calculate the Euler characteristics of the local systems
$\VVl^{\prime}$ on $\md$ if we have the strata and the Euler numbers 
of the strata on $M_3$. It seems very
difficult to determine these Euler numbers directly, but by using
our moduli spaces over finite fields we shall
determine sufficiently many linear relations between them.

\end{section}
\begin{section}{The stratification of the moduli of non-hyperelliptic curves}
\label{sec-nonh}
Also in this section we write $\md$ for ${\mathcal M}_3 \otimes \CC$ and
$\VVl'$ for ${}_0 \VVl'$.
The calculation of the Euler characteristic of $\VVl^{\prime}$ on $\md$
is done by using a stratification on $\md$ (and on $M_3$). First we have
$\md= \md^0 \sqcup \hy$ with $\hy$ the
hyperelliptic locus and $\md^0$ the locus of non-hyperelliptic curves
of genus $3$. Since an algorithm for the calculation of the
Euler characteristic of $\VVl^{\prime}$ on $\hy$ was given
in \cite{BvdG} we are left with calculating $e_c(\md^0,\VVl^{\prime})$.
For this we use the stratification of $\md^0$ given by the
automorphism group in the manner described above.

A non-hyperelliptic curve $C$ of genus $3$ over the field of
complex numbers $\CC$ can be given as a smooth projective curve of
degree $4$ in the projective plane, the canonical image of $C$. 
The automorphism group acts on the vector space $H^0(C,\Omega_C^1)$ 
and thus induces and is in fact induced by a projective automorphism of the 
projective space of $H^0(C,\Omega_C^1)$. It is well-known which 
groups occur as an automorphism group of a non-hyperelliptic 
curve of genus $3$. We list the possibilities in Table \ref{tab1}.
This list is found in \cite{Vermeulen}, cf. also \cite{Magaard}. 

\begin{table}[htbp]\caption{The automorphism groups} \label{tab1} 
$
\begin{array}{|c|c||r|l|c|}
\hline
i & G & \# G & \text{Normal form: } f & \dim(\Sigma(G)) \\
\hline \hline
0 & 1 & 1& f & 6 \\
1 & \ZZ/2\ZZ & 2 & x^4+x^2f(x,y)+g(y,z) & 4 \\
2 & V_4 & 4 & x^4+y^4+z^4+ax^2y^2+by^2z^2+cx^2z^2 & 3 \\
3 & \ZZ/3\ZZ & 3 & yz^3+x(x-y)(x-ay)(x-by) & 2 \\
4 & \mathbb{S}_3 & 6 & x^3y+y^3z+x^2y^2+axyz^2+bz^4 & 2 \\
5 & D_4 & 8 & x^4+y^4+z^4+ax^2y^2+bxyz^2 & 2 \\
6 & \ZZ/6\ZZ & 6 & z^3y+x^4+ax^2y^2+y^4 & 1 \\
7 & \Gamma_{16} & 16 & x^4+y^4+z^4+ax^2y^2 & 1 \\
8 & \mathbb{S}_4 & 24 & x^4+y^4+z^4+a(x^2y^2+x^2z^2+y^2z^2) & 1 \\
9 & \ZZ/9\ZZ & 9 & x^4+xy^3+yz^3 & 0 \\
10 & \Gamma_{48} & 48 & x^4+y^4+z^4+(4\zeta_3+2)x^2y^2 & 0 \\
11 & \Gamma_{96} & 96 & x^4+y^4+z^4 & 0 \\
12 & \Gamma_{168} & 168 & x^3y+y^3z+z^3x & 0 \\
\hline
\end{array}
$
\end{table}

Accordingly, the moduli space $\md^0$ carries a stratification indexed by the
automorphism groups of non-hyperelliptic curves of genus $3$. There are
thirteen (open) strata $\Sigma_{G_i}$ and the inclusion relations between 
the closures of the strata are given in Table \ref{tab2}. 

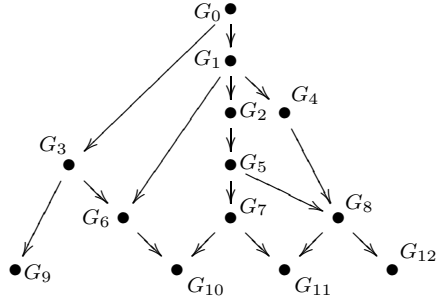
\begin{table}[htbp]\caption{The diagram of strata}\label{diagram} \label{tab2}
$\xymatrix@R=0.7pc@C=0.7pc{&&&&\bullet\ar[d]_<<<<<{G_1}_>>>>>>>>>>>>>>>{G_0}
    \ar[dddlll]_<<<<<<<<<<<<<<<<<<<<<<<<<<<<<<<<<<<<<<<<<<<<<<{G_3} \\
&&&&\bullet\ar[d]^<<<<<<<<<{G_2}\ar[dr]^<<<<<<<{G_4}
    \ar[dddll]^<<<<<<<<<<<<<<<<<<<<<<<<<<<<<<<<<<<<<<{G_6}\\
&&&&\bullet\ar[d]^<<<<<<<<<{G_5}&\bullet
    \ar[ddr]^<<<<<<<<<<<<{G_8}\\
&\bullet\ar[dr]\ar[ddl]_<<<<<<<<<<<<<<<<<<<<<<<{G_9}&&&\bullet\ar[drr]
    \ar[d]^<<<<{G_7}\\
&&\bullet\ar[dr]&&\bullet\ar[dl]^<<<<<<<<<<<<<<<<<<<{G_{10}}
    \ar[dr]&&\bullet\ar[dl]^<<<<<<<<<<<<<{G_{11}}
    \ar[dr]^<<<<<<<<<<<<<{G_{12}}\\
\bullet&&&\bullet&&\bullet&&\bullet}$
\end{table} 

The groups $G_i$ are listed in \cite{Vermeulen} as subgroups of 
$\mathrm{PGL}_3$ acting on the coordinates $x$, $y$, $z$ of the normal forms 
found in Table \ref{tab1}.
In view of Prop.\ \ref{diminv} we need the representations $\rho(G_i)$ on 
the space of differentials $H^0(C,\Omega_c^1)$,
or more precisely, the eigenvalues of all the group elements. 
These representations are listed in Table \ref{tab3}, where the standard basis 
is identified (on $z \neq 0$) with $\omega$, $y_1 \, \omega$, $y_2 \, \omega$, 
where $y_1:=x/z$, $y_2:=y/z$, $G(y_1,y_2):=f(y_1,y_2,1)$ and $\omega$ is the 
differential found by gluing $dy_1/(\partial G / \partial y_2)$ and 
$dy_2/(\partial G / \partial y_1)$. 

The notation for Table \ref{tab1} and \ref{tab3} is as follows: 
$V_4$ is the Klein four-group; 
$D_4$ is the dihedral group of order $8$; 
$\mathbb{S}_n$ denotes the symmetric group on $n$ letters; 
$\Gamma_{n}$ denotes a group of order $n$ 
with $\Gamma_{168}={\rm SL}(3,\FF_2)$;
$\zeta_n$ is a $n$-th root of unity. 

\begin{table}[htbp]\caption{The eigenvalues} \label{tab3} 
$
\begin{array}{|c|r|c|}
\hline
i & \# G & \text{Generators of $\rho(G)$} \\
\hline \hline
0  & 1 & 1 \\
1  & 2 & {\rm diag}(-1,1,-1) \\
2  & 4 & {\rm diag}(-1,1,-1), \,  {\rm diag}(-1,-1,1)\\
3  & 3 & {\rm diag}(\zeta_3^2,\zeta_3,\zeta_3) \\
4  & 6 & \left(\begin{matrix} -1&0&0\\ 0&0&-1 \\ 0&-1&0\\ \end{matrix} \right),\,
{\rm diag}(1,\zeta_3,\zeta_3^2) \\
5  & 8 & \left(\begin{matrix} -1&0&0\\ 0&0&-1\\ 0&-1&0    \\ \end{matrix} \right), \,
{\rm diag}(1,i,-i) \\
6  & 6 & {\rm diag}(-\zeta_3^2,\zeta_3,-\zeta_3) \\
7  & 16 & {\rm diag}(-1,1,-1), \, {\rm diag}(1,i,-i), \,
\left(\begin{matrix} 1&0&0\\ 0&0&-1\\ 0&1&0    \\ \end{matrix} \right) \\
8  & 24 & \left(\begin{matrix} 0&0&1\\ 1&0&0 \\ 0&1&0\\ \end{matrix} \right),\,
\left(\begin{matrix} 1&0&0\\ 0&0&-1 \\ 0&1&0\\ \end{matrix} \right)\\
9  & 9 & {\rm diag}(\zeta_9^2,\zeta_9^4,\zeta_9) \\
10  & 48 & \frac{1}{\sqrt 2} \left(\begin{matrix} \sqrt{2}\,\zeta_3&0&0\\ 0&\zeta_8 &\zeta_8^3 \\
0&\zeta_8 & \zeta_8^7 \\ \end{matrix}\right), \,
\frac{1}{\sqrt 2} \left(\begin{matrix} -\sqrt{2}\,\zeta_3^2 &0&0\\ 0&\zeta_8^5 &\zeta_8 \\ 0&\zeta_8^7 &\zeta_8^7
\\ \end{matrix} \right)\\
11  & 96 & \left(\begin{matrix} 0&0&1\\ 1&0&0 \\ 0&1&0\\ \end{matrix} \right),\,\left(\begin{matrix} 0&0&-i\\ 0&i&0 \\ -1&0&0\\ \end{matrix} \right)\\
12 & 168 & {\rm diag}(\zeta_7, \zeta_7^4, \zeta_7^2), \,
\left((\zeta_7^{2ij}-\zeta_7^{-2ij})/(-\sqrt{-7}) \right)\\
\hline
\end{array}
$
\end{table}

\end{section}
\begin{section}{Linear relations between the Euler numbers}
In this section we use information 
obtained by counting points over finite fields
to get linear relations between the Euler numbers of the strata. 

We consider the moduli space $\md^0$ over $\ZZ[1/2,1/\ell]$.
By the work of the first author the Euler characteristic 
$a_{\lambda}=e_c(\md^0 \otimes \bar{\FF}_p, \VVl^{\prime})$ 
are known for low weight local systems. 
More specifically, in \cite{B1}, for $|\lambda| \leq 7$, the trace 
of Frobenius on $e_c(\md^0 \otimes \bar{k}, \VVl^{\prime})$ is computed 
for all finite fields $k$, and these traces turn out to be given by a
polynomial in the cardinality of the field. 
By Deligne's proof of the Weil conjectures (see \cite{WeilII}) this 
implies that the Euler characteristic is obtained by substituting $1$ 
for the cardinality of the field in this polynomial. 
The results for the traces, in the cases $|\lambda| = 6,7$, are presented 
in Theorems 16.2 and 16.3 of \cite{B1}. For the results in terms of Euler 
characteristics see Table \ref{tab4} and Table \ref{tab5}. 

Since the Euler characteristics 
$e_c(\md^0 \otimes k, \VVl^{\prime})$
for any field $k$ of characteristic different from $2$
are the same
this leads, by \eqref{eq-eulstrat}, to a system of linear equations in the 
unknowns $e_i:=e_c(\Sigma_{G_i})$
$$
\sum_{i=0}^{12} \dim(\VVl^{G_i}) \, e_i = a_{\lambda}.
$$
We shall use these equations for all $\lambda$ with $|\lambda| \leq 4$ 
and $\lambda= (6,0,0)$. We obtain in this way the following result.

\begin{proposition}\label{eulernumbers}
The Euler numbers $e_i=e(\Sigma_{i})$ 
of the strata $\Sigma_{i}$ on $\md^0\otimes \CC$ 
satisfy the following equations:
$$
e_1=-3e_0, \, e_2=2e_0+e_8+1, \, e_4=-e_8-1, \, e_5=-e_8,
$$
and
$$
e_3=0,  \, e_6=e_7=-1, \, e_9=e_{10}=e_{11}=e_{12}=1.
$$
\end{proposition}
Note that this proposition 
agrees with the fact that each of the strata $\Sigma_i$
consists of one point for $i=9,10,11$ and $12$.
These equations do not determine the Euler numbers $e_i$,
but there are some easy relations between the linear equations
in the Euler numbers $e_i$ that reduce the number of unknowns. 
If we write 
$$
k_i(\lambda):=\dim(\VVl^{G_i})
$$
and substitute the
relations of Proposition \ref{eulernumbers} in
$\sum_{i=0}^{12} k_i(\lambda)e_i$, then the coefficients
of $e_0$ and $e_8$ are 
$c_0=k_0(\lambda)-3 k_1(\lambda)+2 k_2(\lambda)$ and
$c_8=k_2(\lambda)-k_4(\lambda)-k_5(\lambda)+k_8(\lambda)$.
The following lemma shows that $c_0=0$ and $-c_0+2c_8=0$,
i.e., the sum $\sum_{i=0}^{12} k_i(\lambda)e_i$ 
does not depend upon $e_0$ and $e_8$.
\begin{lemma}
Let $k_i(\lambda)=\dim(\VVl^{G_i})$. Then we have the relations
$$
k_0(\lambda)-3\, k_1(\lambda)+2\, k_2(\lambda)=0,
$$
and
$$
-k_0(\lambda)+3\, k_1(\lambda)-2\, k_4(\lambda)-2\, k_5(\lambda)+2\, 
k_8(\lambda)=0.
$$
\end{lemma}
\proof 
We use the formula 
$$
\dim(\VVl^G)={\frac{1}{\#G}}  \sum_{\gamma\in G} J_{\lambda}(a_{\gamma},b_{\gamma},c_{\gamma}), 
$$
where 
$\{a_{\gamma},b_{\gamma},c_{\gamma}, 
a_{\gamma}^{-1},b_{\gamma}^{-1},c_{\gamma}^{-1}\}$
are the eigenvalues of $\gamma$ on the standard representation. 
The lemma now follows from a glance at the table of eigenvalues 
from which it follows that the weighted sum of eigenvalues agree. 
The first relation is, for example, clear from the fact that the
eigenvalues for the groups $G_0, G_1$ and $G_2$ are 
$\{[1,1,1]\}$, $\{[1,1,1],[-1,-1,1]\}$ and
$\{[1,1,1],[-1,-1,1],[-1,-1,1],[-1,-1,1]\}$.
\qed 

Combining 
the preceding lemma and proposition gives the following result.

\begin{proposition} For $i=0,\ldots,12$ we let
$k_i(\lambda):=\dim(\VVl^{G_i})$.
Then the Euler characteristic of $\VVl^{\prime}$ on the non-hyperelliptic locus
$\md^0$ is given by
$$e_c(\md^0,\VVl^{\prime}) = k_2(\lambda)- k_4(\lambda)- k_6(\lambda)- 
k_7(\lambda)+ \sum_{j=9}^{12} k_j(\lambda).
$$
\end{proposition}
By the additivity of the Euler characteristic we get the
Euler characteristic on $\md$ by adding the contribution from the 
hyperelliptic locus. An algorithm for determining the
values $e_c({\mathcal H}_3,\VVl^{\prime})$ was given in \cite{BvdG}.

\begin{remark}
Note that for local systems $\VVl^{\prime}$ of odd weight $|\lambda|$
the Euler characteristic vanishes on the hyperelliptic locus
$\hy$ in view of the presence of the hyperelliptic
involution. We therefore present results on the Euler characteristic
of $\md$ for local systems of even and odd weight separately in
Table \ref{tab4} and Table \ref{tab5}. 
Note that the values we obtain for the Euler characteristics on the
hyperelliptic locus (and its complement) are only valid if the 
characteristic is either $0$ or greater than $2$. The closure of 
$\mathcal{H}_3$ in $\overline{\mathcal{M}}_3$ is namely singular
in characteristic $2$, along the closure of the strata consisting of 
a rational backbone with three elliptic tails,
c.f. \cite[Theorem 2.8]{Bertin-Maugeais}. 
For examples where the Euler characteristic 
behaves differently in characteristic $2$, see \cite{B1} or \cite{B2}).
\end{remark}

\begin{subsection}{The hyperelliptic case and the case $g=2$}
It is interesting
to apply the method used here for the non-hyperelliptic case
to the hyperelliptic case. We have a stratification ${\mathcal H}_3$
by strata $\Sigma_{H_i}$ parametrized by the possibilities
$H_i$, $i=1,\ldots,11$, for the automorphism group of a hyperelliptic curve,
cf., \cite{BvdG}. 
The values $b_{\lambda}=e_c({\mathcal H}_3 \otimes \bar{\FF}_p,\VVl^{\prime})$
for $|\lambda| \leq 6$, obtained by counting points over 
finite fields (see \cite{B2}) give us a system of linear equations 
for the Euler
numbers $g_i=e(\Sigma_{H_i})$ with a $3$-dimensional solution space.
Just as for the genus $3$ non-hyperelliptic case, the Euler characteristics
do not depend on these remaining three parameters. 

If the same procedure is applied to ${\mathcal M}_2$ you actually find 
the Euler characteristics of all individual strata from the values 
$c_{\lambda}=e_c({\mathcal H}_2 \otimes \bar{\FF}_p,\VVl^{\prime})$ 
for $|\lambda| \leq 6$. 
\end{subsection}
\end{section}
\begin{section}{Branching}
In order to calculate the Euler characteristic on ${\mathcal A}_3$
we need branching formulas for representations of the symplectic group.
Since the Euler characteristic of $\VVl$ is unchanged by twisting
with the multiplier representation, we may and do work with 
${\rm Sp}(6):={\rm Sp}(6,\QQ)$ instead of ${\rm GSp}(6,\QQ)$.

Let $U$ be a finite-dimensional irreducible complex
representation of the group $G:={\rm Sp}(2)^3\rtimes \mathbb{S}_3$,
where $\mathbb{S}_3$ denotes the symmetric group of the three factors
${\rm Sp}(2):={\rm SL}(2,\QQ) $.
We will need branching formulas from ${\rm Sp}(6)$ to $G$. In the 
literature we found formulas for branching from ${\rm Sp}(6)$ to 
${\rm Sp}(2)^3$ (see, e.g., \cite{lep}), but none to $G$.
The restriction of $U$ to ${\rm Sp}(2)^3$ decomposes as a direct
sum $U=\oplus_{i \in I} W_i$ of irreducible representations of ${\rm Sp}(2)^3$.
The group $\mathbb{S}_3$ acts transitively on the index set. 
Let $\Sigma=\{ \sigma \in \mathbb{S}_3: \sigma(W)=W\}$
be the stabilizer of a component $W=W_i$. Then $W$ is an irreducible 
representation of the semi-direct product ${\rm Sp}(2)^3\rtimes \Sigma$,
hence of the form 
$V_{\alpha} \boxtimes V_{\beta} \boxtimes V_{\gamma}$ with a
$\Sigma$-action, where $V_{n}$ for $n\geq 0$ denotes the irreducible
representation of rank $n+1$ of ${\rm Sp}(2)$. If $\Sigma$
contains a $3$-cycle then $\alpha=\beta=\gamma$ and if $\Sigma$
contains a $2$-cycle then $\# \{ \alpha,\beta,\gamma\}\leq 2$.

If $\Sigma$ is trivial then $\# I =6$ and $U$ consists of six
copies of $V_{\alpha} \boxtimes V_{\beta} \boxtimes V_{\gamma}$.
In this case $U$ is obtained by inducing the
representation $V_{\alpha} \boxtimes V_{\beta} \boxtimes V_{\gamma}$
from ${\rm Sp}(2)^3$ to $G$. We denote it by $R_{\alpha,\beta,\gamma}$.
It is independent of the ordering of $\alpha,\beta,\gamma$ and
therefore we shall assume that $\alpha \geq \beta \geq \gamma$.

Next, suppose that $\Sigma$ is of order $2$. Then $W$ is of the form
$V_{\alpha} \boxtimes V_{\beta}  \boxtimes V_{\beta}$ and $\Sigma$
acts on it, making it a representation
of ${\rm Sp}(2)^3\rtimes \Sigma$. There are two possibilities
for the action of a generator of $\Sigma$:
$$
u \boxtimes v \boxtimes w \mapsto u \boxtimes w \boxtimes v 
\quad \text{\rm or} \quad  -u \boxtimes w \boxtimes v
$$
and $U$ is then the induced representation from ${\rm Sp}(2)^3\rtimes \Sigma$
to $G$. We denote the two possibilities by $R_{\alpha,\beta}^{+}$
and  $R_{\alpha,\beta}^{-}$. Note that $\alpha$ and $\beta$ need not be
different.

If  $\Sigma$ is of order $3$ then  $\alpha=\beta=\gamma$
and a generator $\sigma$ of $\Sigma$ acts on $V_{\alpha} 
\boxtimes V_{\alpha} \boxtimes V_{\alpha}$ via
$$
u \boxtimes v \boxtimes w\mapsto \rho^{\epsilon} \,
 w\boxtimes u \boxtimes v,
$$
with $\rho$ a primitive third root of $1$ and $\epsilon= 0, 1$ or $2$.
The representation $U$ is then the induced representation from 
${\rm Sp}(2)^3\rtimes \Sigma$ to $G$. If $\epsilon = 0$ we denote it by 
$T_{\alpha}$, while for $\epsilon=1$ and $2$ the result is the same and 
is denoted by $T_{\alpha}^{\prime}$.

Finally, if $\Sigma=\mathbb{S}_3$ then $U$ is 
$V_{\alpha} \boxtimes V_{\alpha} \boxtimes V_{\alpha}$  
with the action of $\mathbb{S}_3$ given by 
\begin{multline*}
u_1 \boxtimes u_2 \boxtimes u_3 \mapsto 
u_{\sigma^{-1}(1)} \boxtimes u_{\sigma^{-1}(2)}  \boxtimes u_{\sigma^{-1}(3)}
\quad \text{\rm or} \\
\quad {\rm sgn}(\sigma)
u_{\sigma^{-1}(1)} \boxtimes u_{\sigma^{-1}(2)}  \boxtimes u_{\sigma^{-1}(3)}.
\end{multline*}
We denote the two representations by $R_{\alpha}^{+}$ and $R_{\alpha}^{-}$. 

We have the following two relations:
$$
T_{\alpha}= R_{\alpha}^{+} \oplus R_{\alpha}^{-}
\qquad \text{\rm and} \qquad R_{\alpha,\alpha}^{+} = T_{\alpha}^{\prime}
\oplus R_{\alpha}^{+}.
$$
Thus we obtain the following result.

\begin{lemma} \label{branch}
Every irreducible representation of ${\rm Sp}(2)^3\rtimes \mathbb{S}_3$
is a virtual sum of representations of the form $R_{\alpha,\beta,\gamma}$,
$R_{\alpha,\beta}^{\pm}$ and $R_{\alpha}^{\pm}$.
\end{lemma}

The exterior products $\wedge^i V_{1,0,0}$ for $i=1,2,3$ of the standard
representation of ${\rm Sp}(6)$ form a basis for the representation
ring of ${\rm Sp}(6)$. Their restrictions to the subgroup 
${\rm Sp}(2)^3\rtimes \mathbb{S}_3$ are given by
$$
V_{1,0,0}|G=R_{1,0}, \quad
\wedge^2 V_{1,0,0}|G= R_{0,1}^{-} \oplus R_{0,0}^{+}, \quad
\wedge^3 V_{1,0,0}|G=R_{1,1}^{+} \oplus R_{1,0,0}.
$$
Explicitly, the representation $V_{\lambda}$ corresponds to the symmetric
function $J_{\lambda}$ and by writing $J_{\lambda}$ in terms of the
elementary symmetric functions $e_1$, $e_2$ and $e_3$ one finds the
expression of $V_{\lambda}$ in terms of the three exterior products.
Thus we can determine the branching completely. We implemented this
in a computer program.
\end{section}
\begin{section}{The Euler characteristics on $\A_3$} \label{eulA3}
In this section we will use the stratification, 
$$\A_3 = t_3(\md^0) \; \sqcup \; t_3(\mathcal{H}_3) \; 
\sqcup \; t_2(\mathcal{M}_2) \times \A_1 \;\sqcup \;\A_{1,1,1},
$$
to compute 
$e_c(\A_3,\VVl)$ for any partition $\lambda$, where 
$t_g\colon \mathcal{M}_g \to \A_g$ is the Torelli morphism. 
For hyperelliptic curves the Torelli morphism has degree one 
and hence $e_c(t_3(\mathcal{H}_3),\VVl)$ 
(resp.\ $e_c(t_2(\mathcal{M}_2),\VVl)$)
is equal to $e_c(\mathcal{H}_3,\VVl')$ (resp.\
$e_c(\mathcal{M}_2,\VVl')$). On the locus $\md^0$ in 
$\mathcal{M}_3$ of non-hyperelliptic curves the Torelli morphism 
has degree two, due to the existence of the automorphism ``$-1$'' 
on the Jacobians. That is, if $G$ is the automorphism group of a 
non-hyperelliptic curve $C$, then $G':=\{\pm g:g \in G\}$ is the 
automorphism group of its Jacobian. 
The canonical isomorphism between $H^1(C)$ and $H^1(J(C))$ and the 
fact that $-1 \in G'$ acts as $-1$ on $H^1(J(C))$ then shows how the 
elements of $G'$ act on $H^1(J(C))$. The Torelli morphism is a homeomorphism 
of the coarse moduli spaces
and thus the Euler characteristic of the strata of $t_3(\md^0)$ 
are the same as the corresponding strata of $\md^0$. 
Since the polynomials $J_{\lambda}$ are even (resp.\ odd) if 
$\lambda$ is even (resp.\ odd), we find by applying formula 
\eqref{eq-eulstrat}, that $e_c(t_3(\md^0),\VVl)$ is equal to 
$e_c(\md^0,\VVl')$ if $\lambda$ is even and zero if $\lambda$ is odd. 

We will apply a K\"unneth formula to compute 
$e_c(t_2(\mathcal{M}_2) \times \A_1,\VVl)$. 
For this we need to know the branching, that is the restriction 
of the representation $V_{\lambda}$ of $\Sp$ to 
${\rm Sp}_{4,2}:=\mathrm{Sp}(4,\QQ) \times \mathrm{Sp}(2,\QQ)$. 
If $V_{\lambda}|_{{\rm Sp}_{4,2}} =\bigoplus_{\mu,\nu} (V_{\mu} 
\boxtimes V_{\nu})^{\oplus m_{\mu,\nu}}$ it follows that 
$$
e_c(t_2(\mathcal{M}_2) \times \A_1,\VVl)= \sum_{\mu,\nu} m_{\mu,\nu} 
\, e_c(\mathcal{M}_2,t^*\VV_{\mu}) \, e_c(\A_1,\VV_{\nu}). 
$$
A formula for the branching from $\Sp$ to ${\rm Sp}_{4,2}$ 
(i.e., for the numbers $m_{\mu,\nu}$) can for instance be found in 
\cite{lep}.

Finally, to compute $e_c(\A_{1,1,1},\VVl)$ we will also use branching, 
in this case from $\Sp$ to $G=\mathrm{Sp}(2,\QQ)^3 \rtimes \mathbb{S}_3$, 
where the symmetric group $\mathbb{S}_3$
acts by permuting the three factors ${\rm Sp}(2,\QQ)$. This branching was 
treated in the preceding section. 
Since $\A_{1,1,1} = (\A_1)^{3} / \mathbb{S}_3$, we have to
compute the invariant part of the cohomology of the local system 
$\VVl$ restricted to $(\A_1)^{3}$. 

By Lemma \ref{branch} we can write $V_{\lambda}|G$, as a virtual 
representation,
$$
\sum_{\alpha \geq \beta \geq \gamma} m_{\alpha,\beta, \gamma}
\, R_{\alpha,\beta, \gamma}+ 
\sum_{\alpha,\beta} \left( m_{\alpha,\beta}^{+} \, R_{\alpha,\beta}^{+}
+m_{\alpha,\beta}^{-} \, R_{\alpha,\beta}^{-}\right) +\sum_{\alpha} \left( m_{\alpha}^{+} \, R_{\alpha}^{+} +
m_{\alpha}^{-} \, R_{\alpha}^{-}\right).
$$ 
We then have
\begin{multline*}
e_c(\A_{1,1,1},\VVl)= 
\sum_{\alpha \geq \beta \geq \gamma} m_{\alpha,\beta,\gamma} 
\, e_c(\A_1,\VV_{\alpha}) \, e_c(\A_1,\VV_{\beta})
\, e_c(\A_1,\VV_{\gamma}) \\
+\sum_{\alpha ,\beta} \Bigl( m_{\alpha,\beta}^{+} 
\, \binom{e_c(\A_1,\VV_{\beta})+1}{2} + 
m_{\alpha,\beta}^{-} \, \binom{e_c(\A_1,\VV_{\beta})}{2}\Bigr)
 \, e_c(\A_1,\VV_{\alpha})\\
+\sum_{\alpha} \Bigl( m_{\alpha}^{+} \, \binom{e_c(\A_1,\VV_{\alpha}) +2}{3}+ m_{\alpha}^{-} \, \binom{e_c(\A_1,\VV_{\alpha})}{3} \Bigr).  
\end{multline*}

Thus our knowledge of the Euler characteristics of local systems
for $g=1$ and $g=2$ together with those of local systems on $\md$
suffices to calculate algorithmically
the Euler characteristics on ${\mathcal A}_3$.

We calculated these for all $\lambda$ with $|\lambda|\leq 60$, 
see Table \ref{tab7} for the results for $|\lambda|\leq 18$. 
A first check is that the value $e_c(\A_3,\VVl)=5$ for
the trivial local system given by $\lambda=(0,0,0)$ 
agrees with a result of Hain, who
calculated the rational cohomology of $\A_3$, cf., \cite{Hain}. 
A further indication
of the correctness, besides the fact that while we are summing
rational numbers we always found integer values,
is that the absolute value of the
Euler characteristic of $\VVl$ on $\md$ is in general much smaller than
the Euler characteristic of $\VVl$ on the hyperelliptic locus and its
complement. This is illustrated in Table \ref{tab4}.
A similar phenomenon was observed for the genus $g=2$ case (cf., \cite{FvdG})
and can be observed for $\A_3$ too, as illustrated in Table \ref{tab6}.
\end{section}
\begin{table} \caption{Euler characteristics of $\VVl^{\prime}$ 
(even weight) on
the parts of the stratification of $\md$.} \label{tab4}
$
\begin{array}{|c|r|r|r||c|r|r|r|}
\hline \lambda & \hy & \md^0 & \md & \lambda & \hy & \md^0 &
\md \\
\hline \hline
(0,0,0) & 1  & 2 & 3 & (5,2,1) & -10 & 6 & -4 \\
(2,0,0) & -1 & 1 & 0 & (4,4,0) & -5 & 7  & 2 \\
(1,1,0) & 0&-1 & -1 & (4,3,1) & -4 & 4 & 0\\
(4,0,0) & -1& 1 & 0  & (4,2,2) & -7 & 8 & 1\\
(3,1,0) & 0 & 0& 0 & (3,3,2) & -2 & 0 & -2\\
(2,2,0) & -1& 1 & 0  &(10,0,0) & -17 & 13 & -4 \\
(2,1,1) & 0& 1& 1 & (9,1,0) & -22& 20& -2\\
(6,0,0) & -5& 4& -1 &(8,2,0) & -43 & 37& -6\\
(5,1,0) & -2 & 1 & -1 & (8,1,1) & -8 & 7 &-1\\
(4,2,0) & -5 & 3 & -2 & (7,3,0) & -34 & 26 & -8\\
(4,1,1) & 0 &-1 & -1 & (7,2,1) & -32 & 28 & -4\\
(3,3,0) & 0 & 2 & 2 & (6,4,0) & -37& 31& -6\\
(3,2,1) & 0& 0 & 0 & (6,3,1) & -26 & 17 & -9\\
(2,2,2) & -3 & 4 & 1 & (6,2,2) & -27 & 25 & -2\\
(8,0,0) & -7 & 9 & 2 & (5,5,0)& -6 & 8 & 2\\
(7,1,0) & -8 & 7 & -1 & (5,4,1)& -22 & 22 & 0\\
(6,2,0) & -13 & 14 & 1 & (5,3,2)& -12 & 9 & -3\\
(6,1,1) & -2 & 0 & -2 & (4,4,2)& -15 & 13 & -2\\
(5,3,0) & -10 & 8 & -2 &(4,3,3) &  0& -2 & -2 \\
\hline
\end{array}
$
\end{table}

\begin{table} \caption{Euler characteristics of $\VVl^{\prime}$ 
(odd weight) on $\md$.} \label{tab5}
$
\begin{array}{|c|r||c|r||c|r||c|r|}
\hline
\lambda & \md & \lambda & \md & \lambda & \md & \lambda & \md\\
\hline \hline
(1,0,0) & 0  & (2,2,1) & 0  & (3,2,2) & 0  & (5,3,1) & 20 \\
(3,0,0) & 0  & (7,0,0) & -2 & (9,0,0) & 4  & (5,2,2) & 2 \\
(2,1,0) & 0  & (6,1,0) & 4  & (8,1,0) & 8  & (4,4,1) & 2 \\
(1,1,1) & 0  & (5,2,0) & 4  & (7,2,0) & 10 & (4,3,2) & 4\\
(5,0,0) & 0  & (5,1,1) & 10 & (7,1,1) & 18 & (3,3,3) & 8 \\
(4,1,0) & 2  & (4,3,0) & 2  & (6,3,0) & 20 & (11,0,0)& 4 \\
(3,2,0) & -2 & (4,2,1) & 2  & (6,2,1) & 12 & (10,1,0) & 30 \\
(3,1,1) & 2  & (3,3,1) & 4  & (5,4,0) & 6  & (9,2,0) & 36 \\
\hline
\end{array}
$
\end{table}

\begin{table} \caption{Euler characteristics of $\VVl$ (of high weight) on
the parts of the stratification of $\A_3$.} \label{tab6}
$
\begin{array}{|c|r|r|r|r|r|}
\hline
\lambda & {\mathcal H}_3 & {\md^0} & {\mathcal M}_2 \times A_1&
\A_{1,1,1} & \A_{3}\\
\hline\hline
(40,0,0) & -3825 & 3257 &  731 & -161 & 2\\
(32,5,3) & -188587 & 156651  & 44843 & -12615 & 292\\
(24,12,4) & -502970 & 419733 & 116703& -32415 & 1051\\
(21,15,4) & -351374 & 292508 & 82372& -22910 & 596\\
(14,13,13) & -2262 &  1795 & 649& -191 & -9\\
\hline
\end{array}
$
\end{table}

\begin{table} \caption{Euler characteristics of $\VVl$ on $\A_3$.}
\label{tab7}
$
\begin{array}{|c|r||c|r||c|r||c|r|}
\hline
\lambda & \A_{3} & \lambda & \A_{3} & \lambda & \A_{3} & \lambda & \A_{3}\\
\hline\hline
(0,0,0) & 5   &  (4,3,3) & 0   &  (7,6,1) & 0   &  (18,0,0) & -5\\
(2,0,0) & -2  &  (12,0,0) & -1 &  (7,5,2) & 1   &  (17,1,0) & -7\\
(1,1,0) & -1  &  (11,1,0) & -3 &  (7,4,3) & 0   &  (16,2,0) & -4\\
(4,0,0) & -2  &  (10,2,0) & 0  &  (6,6,2) & 4   &  (16,1,1) & -1 \\
(3,1,0) & 0   &  (10,1,1) & -1 &  (6,5,3) & 1   &  (15,3,0) & 0\\
(2,2,0) & 0   &  (9,3,0) & 0   &  (6,4,4) & 0   &  (15,2,1) & 0\\
(2,1,1) & 1   &  (9,2,1) & 0   &  (5,5,4) & 0   &  (14,4,0) & -4\\
(6,0,0) & -3  &  (8,4,0) & 0   &  (16,0,0) & -2 &  (14,3,1) & -4\\
(5,1,0) & -1  &  (8,3,1) & -2  &  (15,1,0) & -4   &  (14,2,2) & 7\\
(4,2,0) & 0   &  (8,2,2) & 3   &  (14,2,0) & 4    &  (13,5,0) & 0\\
(4,1,1) & 1   &  (7,5,0) & -2  &  (14,1,1) & -1 &  (13,4,1) & 0\\
(3,3,0) & 0   &  (7,4,1) & 0   &  (13,3,0) & -4 &  (13,3,2) & 1\\
(3,2,1) & 0   &  (7,3,2) & 1   &  (13,2,1) & 0  &  (12,6,0) & 3\\
(2,2,2) & 1   &  (6,6,0) & 5   &  (12,4,0) & 4  &  (12,5,1) & -3\\
(8,0,0) & 0   &  (6,5,1) & -1  &  (12,3,1) & -3 &  (12,4,2) & 4\\
(7,1,0) & -1  &  (6,4,2) & 0   &  (12,2,2) & 4  &  (12,3,3) & -1\\
(6,2,0) & 1   &  (6,3,3) & 1   &  (11,5,0) & -2 &  (11,7,0) & -3\\
(6,1,1) & 0   &  (5,5,2) & -1  &  (11,4,1) & 0  &  (11,6,1) & 0\\
(5,3,0) & 0   &  (5,4,3) & 0   &  (11,3,2) & -1 &  (11,5,2) & 7\\
(5,2,1) & 0   &  (4,4,4) & 0   &  (10,6,0) & 4  &  (11,4,3) & 0\\
(4,4,0) & 0   &  (14,0,0) & -2 &  (10,5,1) & -2 &  (10,8,0) & 0\\
(4,3,1) & 0   &  (13,1,0) & -3 &  (10,4,2) & 4  &  (10,7,1) & -1 \\
(4,2,2) & 1   &  (12,2,0) & -1 &  (10,3,3) & 0  &  (10,6,2) & 5\\
(3,3,2) & 0   &  (12,1,1) & 0  &  (9,7,0) & 0   &  (10,5,3) & 0\\
(10,0,0) & -4 &  (11,3,0) & -2 &  (9,6,1) & 0   &  (10,4,4) & 4\\
(9,1,0) & 0   &  (11,2,1) & 0  &  (9,5,2) & 0   &  (9,9,0) & 4\\
(8,2,0) & 0   &  (10,4,0) & 4  &  (9,4,3) & 0   &  (9,8,1) & 0\\
(8,1,1) & 1   &  (10,3,1) & -2 &  (8,8,0) & 8   &  (9,7,2) & 0\\
(7,3,0) & -2  &  (10,2,2) & 1  &  (8,7,1) & 1   &  (9,6,3) & 8\\
(7,2,1) & 0   &  (9,5,0) & 0   &  (8,6,2) & 3   &  (9,5,4) & 0\\
(6,4,0) & 0   &  (9,4,1) & 0   &  (8,5,3) & 0   &  (8,8,2) & 3\\
(6,3,1) & -1  &  (9,3,2) & 0   &  (8,4,4) & 8   &  (8,7,3) & -3\\
(6,2,2) & 0   &  (8,6,0) & 0   &  (7,7,2) & 0   &  (8,6,4) & 0\\
(5,5,0) & 0   &  (8,5,1) & -2  &  (7,6,3) & 0   &  (8,5,5) & 0\\
(5,4,1) & 0   &  (8,4,2) & 0   &  (7,5,4) & 0   &  (7,7,4) & -1\\
(5,3,2) & -1  &  (8,3,3) & 0   &  (6,6,4) & 3   &  (7,6,5) & 0\\
(4,4,2) & 0   &  (7,7,0) & -1  &  (6,5,5) & 0   &  (6,6,6) & 3\\
\hline
\end{array}
$
\end{table}

\end{document}